\documentclass{aptpub}

\authornames{Ramtin Pedarsani, Jean Walrand} 
\shorttitle{Stability of Mutliclass Queueing Networks} 

\newcommand{\bW}{\ensuremath{\mathbf{W}}}
\newcommand{\bD}{\ensuremath{\mathbf{D}}}

\newcommand{\bE}{\ensuremath{\mathbf{E}}}
\newcommand{\bM}{\ensuremath{\mathbf{M}}}
\newcommand{\bQ}{\ensuremath{\mathbf{Q}}}

\newcommand{\bR}{\ensuremath{\mathbf{R}}}
\newcommand{\bX}{\ensuremath{\mathbf{X}}}
\newcommand{\bI}{\ensuremath{\mathbf{I}}}
\newcommand{\bt}{\ensuremath{\mathbf{t}}}
\newcommand{\bT}{\ensuremath{\mathbf{T}}}

\newcommand{\be}{\ensuremath{\mathbf{e}}}

\newcommand{\calG}{\ensuremath{\mathcal{G}}}

\newcommand{\calD}{\ensuremath{\mathcal{D}}}
\newcommand{\calJ}{\ensuremath{\mathcal{J}}}

\usepackage{graphicx,epsfig,graphics,epsf,color,subcaption,caption}

\begin{document}

\title{Stability of 
Multiclass Queueing Networks \\under
Longest-Queue and \\ Longest-Dominating-Queue Scheduling} 

\authoronetwosame{Ramtin Pedarsani} 
\authorone[University of California, Berkeley]{Jean Walrand}

\address{Department of Electrical Engineering and Computer Sciences, University of California, Berkeley, CA, 94720, USA}

\emailtwo{ramtin@eecs.berkeley.edu} 
\emailone{wlr@eecs.berkeley.edu}

\begin{abstract}
We consider the stability of robust scheduling policies for multiclass queueing networks. These  are open networks with arbitrary routing matrix and several disjoint groups of queues in which at most one queue can be served at a time. The arrival and potential service processes and routing decisions at the queues are independent, stationary and ergodic. A scheduling policy is called {\em robust} if it does not depend on the arrival and service rates nor on the routing probabilities. A policy is called {\em throughput-optimal} if it makes the system stable whenever the parameters are such that the system can be stable.
We propose two robust polices: longest-queue scheduling and a new policy called longest-dominating-queue scheduling. We show that longest-queue scheduling is throughput-optimal for two groups of two queues. We also prove the throughput-optimality of longest-dominating-queue scheduling when the network topology is acyclic, for an arbitrary number of groups and queues. 
\end{abstract}







\section{Introduction}\label{sec:intro}
We consider the scheduling of multiclass queueing networks \cite{LK91}, \cite{dai99}.  These are queueing networks with disjoint groups of queues that cannot be scheduled simultaneously. MaxWeight scheduling, proposed in \cite{TE92,DL05}, is known to be throughput-optimal for these networks. However, MaxWeight scheduling suffers from high complexity and also dependency on the knowledge of all of the service rates, queue lengths, and routing probabilities. A natural low complexity scheduling algorithm is longest-queue (LQ) scheduling, studied in \cite{TE92a}, \cite{Mc95} and \cite{DW06}, which we discuss in this paper in detail.  

\begin{figure}
\centering
    \includegraphics[width= 0.35\textwidth]{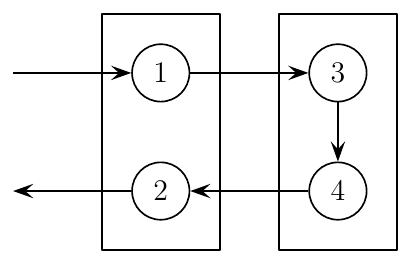}
  \caption{\label{fig0}Lu-Kumar network.} 
\end{figure}

In general, we know that the utilization being less than one for each server is necessary but not sufficient for the stability of a queueing network. This condition specifies that  work arrives at each server at rate less than one. In \cite{LK91}, Lu and Kumar provided an example of a network (Figure \ref{fig0}) with priority scheduling that is unstable despite satisfying the utilization condition. The priority is given to queue 2 in group 1 and to queue 3 in group 2. To see this, assume that $\mu_1 > \mu_3$ and $\mu_4 > \mu_2$ and that queue 3 is initially empty while queue 2 is not.  Group 1 serves queue 2, so that queue 1 is not served and queue 3 remains empty.  Eventually, queue 2 becomes empty and queues 1 and 3 get served until queue 3 becomes empty (because $\mu_1 > \mu_3$).  At that time, queues 4 and 2 get served, until queue 2 becomes empty (because $\mu_4 > \mu_2$).  Thus, queues 2 and 3 are never served together.  Consequently, they form a virtual group and the system cannot be stable unless the utilization of that virtual group is less than one, which is an additional condition not implied by the original utilization condition. The notion of virtual group was introduced in \cite{DV00}. For more details about the proofs, see Lemma 4.1.2 in \cite{dai99}.

If priority is given to the last buffer in each group (queues $2$ and $4$) or first buffer in each group (queues $1$ and $3$), the network is stable. It is shown in  \cite{DW96, KK96, KM95} that, for a re-entrant line, the last-buffer-first-serve and first-buffer-first-serve discipline with pre-emption is throughput-optimal.

We analyze the stability of the queueing networks under robust policies using their fluid model. The fluid model converts a stochastic system into a deterministic system,  based on the functional strong law of large numbers, \cite{malyshev, dai95, stolyar, dai99, Bramson}. Under weak assumptions, the stability of the fluid model implies the stability of the queueing network. To show the stability of the fluid model, we find a Lyapunov function for the differential equations of the system.

There has been relatively little work on the stability of LQ scheduling in the literature. We can  mainly mention the following papers. In \cite{KGL}, Kumar et al. consider sufficient conditions for stability of stable marriage scheduling algorithms in input-queued switches. In \cite{DW06}, Dimakis and Walrand identify new sufficient conditions for longest-queue-first (LQF) scheduling to be throughput-optimal. They combine properties of diffusion-scaled path functionals and local fluid limits into a sharper characterization of stability. See also \cite{LS} for a variation on the first order sufficient condition of that paper (resource pooling). In \cite{BT11}, Baharian and Tezcan consider LQF scheduling for parallel server systems. It is shown that the utilization condition  is sufficient to prove stability if the underlying graph of the parallel server system is a tree. Furthermore, additional drift conditions are provided for the stability of a special parallel server system known as X-model. The network model that we consider is different from all the previous works on LQ scheduling.

In addition to LQ scheduling, we propose a new scheduling policy called longest-dominating-queue (LDQ). This policy is closely related to LQ scheduling. According to this policy, none of the queues that feed a larger maximum-length queue in the network is served.  Among those queues that are not dominated by a maximum-length queue in the network, the longest one is served in each group. We use the fluid model and show that the maximum of the queue lengths is a Lyapunov function to prove the stability of LDQ scheduling in a general acyclic network.

We use boldface letters to denote vectors and matrices. $[.]^T$ denotes matrix transposition operation and $[.]^{-T}$ denotes the matrix inversion and transposition. Let $\mathbf{V}$ and $\mathbf{W}$ be matrices of size $L \times N$, then $\mathbf{V} \preceq \mathbf{W}$ means $v_{ij} \leq w_{ij},$ for all $ i=1, \dots, L$ and $j=1, \dots,N$.  Also,  $\mathbf{1}_L$ and $\mathbf{0}_L$ stand for column vectors of $1$'s and $0$'s with length $L$, respectively, and $\mathbf{0}_{L \times N}$ stands for the  $L \times N$ matrix of all $0$ entries. The indicator function of set $A$ is shown as $1_A$.

The rest of the paper is organized as follows. In Section \ref{sec:prob}, we provide the precise network model and problem formulation. In Section \ref{sec:LQF}, we focus on the LQ scheduling and prove that it is throughput-optimal for multiclass queueing networks with two groups of two queues. In Section \ref{sec:LDQ}, we propose the LDQ scheduling and prove that this new policy is throughput-optimal if the network topology is acyclic.

\section{Problem definition and network model}\label{sec:prob}
In this section, we first introduce the queueing network model that we consider. We discuss the fluid model to analyze the stability of the queueing network. We also provide the utilization condition or stability condition of the network. 

\subsection{Network model}

We consider a  network with $K$ queues and a routing matrix $\bR_{K \times K}$. The entry $r_{ik}$ is the probability that a job goes to queue $k$ upon leaving queue $i$ . Therefore, a job leaves the network upon leaving queue $i$ with probability $ 1 - \sum_k r_{ik}$.  In each queue, jobs are served in their order of arrival. A queueing network is ``acyclic'' if a job cannot visit any queue more than once; that is, if there exists no finite sequence of queues $(i_1,i_2,\ldots,i_\ell,i_1)$ such that $r_{i_1i_2} r_{i_2i_3}\cdots r_{i_\ell i_1} > 0$.  

\begin{figure}
\centering
    \includegraphics[width= 0.4 \textwidth]{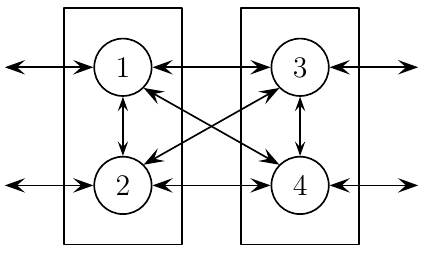}
  \caption{\label{fig1} Multiclass queueing network with two groups of two queues.} 
\end{figure}

The random length of queue $i$ at time $t \geq 0$ is ${\bar X}_i(t)$, the vector of service rates of the queues is $\bmu = [\mu_1, \mu_2, \dots, \mu_K]^T$, and the vector of exogenous arrival rates to the queues is $\blambda = [\lambda_1,\lambda_2,\dots,\lambda_K]^T$. The exogenous arrival processes to the queues are  independent stationary ergodic processes. In this network, not all the queues can be served at the same time. The queues are partitioned into $J$ disjoint groups $\{ \calG_j \}_{j=1}^J$  and only one queue can be served in each group at a time. 
Figure \ref{fig1} illustrates a multiclass queueing network with  $2$ groups of two queues: $\calG_1 = \{1,2 \}$ and $\calG_2 = \{3, 4\}$.

We assume the network is open, i.e., all the jobs eventually leave the network. Since the network is open,
$$\bQ = \bI + \bR^T + (\bR^T)^2 + \dots = (\bI - \bR^T)^{-1} $$
is a finite positive matrix. So, the matrix $(\bI - \bR^T)$ is  invertible. 

The service disciplines that we consider are independent of $\blambda$, $\bmu$ and $\bR$. The goal is to analyze the stability of this network when the utilization condition that will be stated in Section \ref{sec:st} holds. Let ${\bar X}_i(t)$ be the length of queue $i$ at time $t$ and ${\bar T}_i(t)$ be the random total amount of time that queue $i$ has been scheduled up to time $t$. We use the fluid model to analyze the stability of the network. 
Let $X_i(t)$ be the fluid level of queue $i$ at time $t$, and let $T_i(t)$ be the cumulative amount of time that queue $i$ is served up to time $t$ in the fluid model. Then the fluid model equations are 
\begin{align}\label{eq:fluid}
X_k(t) = X_k(0) + \lambda_k t - \mu_k T_k(t)+ \sum_{i=1}^K r_{ik} \mu_{i} T_i(t), \quad k = 1, \dots, K.
\end{align}  
 
The fluid model is {\em stable} if there exits some $\delta >  0$ such that, for each fluid solution with $\|\bX(0)\| \leq 1$, one has $\bX(t)=\mathbf{0}$ for $t>\delta$.  Under weak conditions, the stability of the fluid model implies the stability of the queueing network (e.g., the Harris recurrence of a Markov model in the case of renewal arrival processes, i.i.d. service times and Markov routing).  See \cite{dai95, dai99} for a discussion of such results.
%

\subsection{Utilization condition}\label{sec:st}
Define the drift matrix $\bD$ of the network by
\begin{align}\label{eq:D1}
\bD = \bM (\bR - \bI),
\end{align} where $\bM = \text{diag} \{ \bmu \}$. Since $\bQ  :=(\bI - \bR^T)^{-1}$ is a positive matrix, one has
\begin{equation}\label{eq:D2}
\bD^{-T} = - \bM^{-1} \bQ \preceq \mathbf{0}_{K \times K}.
\end{equation}  
Furthermore, let the vector of nominal traffic of the queues be $\bnu = [\nu_1,\nu_2,\dots,\nu_K]^T$. We have
\begin{equation*}
\nu_i = \lambda_i + \sum_{k=1}^K r_{ki} \nu_k.
\end{equation*}
In matrix form, 
\begin{align}\label{eq:nu}
\bnu = (\bI - \bR^T)^{-1} \blambda.
\end{align}

By Proposition $2.5.3.$ in \cite{dai99}, the necessary condition for the stability of the network (utilization condition) is 
\begin{align}\label{eq:stability}
\sum_{i \in \calG_j} \frac{\nu_i}{\mu_i} \leq 1 - \epsilon, ~ ~ \forall j = 1, \dots J,
\end{align}
for some $\epsilon > 0$.
The utilization condition expresses that work for each group arrives at the network at rate less than one.

Using Equations \eqref{eq:D2} and \eqref{eq:nu}, we find an equivalent stability condition which is ${\mathbf{\tilde{e}}_j}^T\bD^{-T}\blambda + 1 \geq \epsilon $ for all $j = 1, \ldots, J$, where $\mathbf{\tilde{e}}_j$ is a column vector of length $K$ corresponding to group $j$ such that component $i$ of $\mathbf{\tilde{e}}_j$ is
$\tilde{e}_j(i) = 1_{ \{ i \in \calG_j\} }$.

\section{LQ scheduling}\label{sec:LQF}
In this section, 
we define LQ scheduling policy and analyze it in the fluid model. We show that LQ scheduling is throughput-optimal in the queueing network with two groups of two queues.

\subsection{LQ Scheduling for Multiclass Queueing Networks}\label{sec:general}

We provide the fluid model equations for LQ scheduling.  Under this discipline, the longest queue in each group is served. To break ties, we can use a static priority scheduling among the maximum queues, or serve any of the maximum-length queues randomly. 

Let 
$$
S_j(t) = \{i: i \in {\cal G}_j, X_i(t) = \max_{k \in {\cal G}_j} X_k(t)\}
$$ 
be the set of queues with maximum fluid level in group $j$ 
at time $t$. 

Fluid model equations are those defined by \eqref{eq:fluid} together with 
\begin{align}\label{eg:TS}
\sum_{i \in S_j} \dot{T}_i(t) & = 1 \nonumber \\
\sum_{i \in \calG_j \setminus S_j} \dot{T}_i(t) & = 0,
\end{align}
which states that the server in each group spends all its time serving the longest queues.  The main result of this paper is the following theorem.

\begin{thm}\label{thm:1}
The multiclass queueing network with two groups of two queues, $\calG_1 = \{ 1,2 \}$ and $\calG_2 = \{ 3,4 \}$, is stable under LQ scheduling if the utilization condition stated in \eqref{eq:stability} holds.
\end{thm}

The network is illustrated in Figure \ref{fig1}. The rest of this section establishes the proof of Theorem \ref{thm:1}. The key idea for the proof is to find a piece-wise linear Lyapunov function. The stability of multiclass queueing networks under LQ scheduling for general networks remains to be an open problem. 

Let $\bW(t) = \bM^{-1} \bQ \bX(t) = -\bD^{-T} \bX(t)$. Note that $\bQ \bX(t)$ denotes the vector of potential fluid level for each queue at time $t$, and $W_i(t)$ is the potential work to be done at queue $i$ at time $t$. Let 
\begin{align}\label{a}
a = \frac{\mu_1^{-1} (q_{13} + q_{14}) + \mu_2^{-1} (q_{23} + q_{24})}{\mu_3^{-1} (q_{33} + q_{34}) + \mu_4^{-1} (q_{43} + q_{44})},
\end{align}
and
\begin{align}\label{b}
b = \frac{\mu_1^{-1} (q_{11} + q_{12}) + \mu_2^{-1} (q_{21} + q_{22})}{\mu_3^{-1} (q_{31} + q_{32}) + \mu_4^{-1} (q_{41} + q_{42})},
\end{align}
where $\bQ = [q_{ij}]$. 
In Lemma \ref{lem:ab}, we show that $a \leq b$. Let $\beta \in [a,b]$. We show that
\begin{align}\label{V}
V(\bX) = \max(f_1, f_2),
\end{align}
where $f_1 = W_1 + W_2 = {\mathbf{\tilde{e}}_1}^T \bW$ and $f_2 = \beta[W_3 + W_4] = \beta {\mathbf{\tilde{e}}_2}^T \bW$,
is a Lyapunov function. 

\begin{lem}\label{lem:ab}
The following inequality holds:
\begin{align}
a \leq b,
\end{align}
where $a$ is given in \eqref{a} and $b$ is given in \eqref{b}.
\end{lem}

\begin{proof}
It is enough to show the following 4 inequalities:
\begin{align} \label{eq:det1}
(q_{11} + q_{12})(q_{33} + q_{34}) &\geq (q_{13} + q_{14})(q_{31} + q_{32}) \\ \label{eq:det2}
(q_{11} + q_{12})(q_{43} + q_{44}) &\geq (q_{13} + q_{14})(q_{41} + q_{42}) \\
(q_{21} + q_{22})(q_{33} + q_{34}) &\geq (q_{23} + q_{24})(q_{31} + q_{32}) \\ \label{eq:det3}
(q_{21} + q_{22})(q_{43} + q_{44}) &\geq (q_{23} + q_{24})(q_{41} + q_{42}). 
\end{align}
We prove \eqref{eq:det1}. The proof of \eqref{eq:det2}--\eqref{eq:det3} is identical with change of indices. 

Recall that $\bQ = (\bI - \bR^T)^{-1}$. $\bI - \bR$ is an M-matrix \cite{M-matrix}. That is, $\bI - \bR$ has non-negative diagonal entries, non-positive off-diagonal entries, and the sum of entries of each row is non-negative. M-matrices have positive principal minors \cite{M-matrix}; thus, $\det(\bI - \bR)$ is strictly positive, as the network is open. The following equality can be computed: 
\begin{align}
&\det \left( \sqrt{\det(\bI - \bR)} \left [
\begin{array}{cc}
q_{11} + q_{12} & q_{13} + q_{14} \\
q_{31} + q_{32} & q_{33} + q_{34}
\end{array}
\right]
 \right) \nonumber \\
& \qquad \qquad = (1 + r_{21} - r_{22})(1 + r_{43} - r_{44}) + (r_{23} - r_{24})(r_{42} - r_{41}).
\end{align}

It is now sufficient to show that the right-hand side of the above equality is non-negative. This can be shown as follows.
\begin{align*}
&(1 + r_{21} - r_{22})(1 + r_{43} - r_{44}) + (r_{23} - r_{24})(r_{42} - r_{41}) \\
&\qquad \qquad \geq (1 - r_{22})(1  - r_{44}) - (r_{23} + r_{24})(r_{42} + r_{41}) \geq 0.
\end{align*}
The last inequality is due to the facts $1 - r_{22} \geq r_{23} + r_{24}$ and $1 - r_{44} \geq r_{42} + r_{41}$. This completes the proof of the lemma.
\end{proof}

Now we show that $V(\bX)$ is a Lyapunov function. First, it is clear that $V(\bX) = \mathbf{0}$ if and only if $\bX = 0$. We show that at a regular point $t$, $\dot{V}(t) \leq -\delta$ if $\bX(t) \neq \mathbf{0}$ for some $\delta > 0$ by considering 3 different cases. A time $t$ is regular if each component of $\bX(t)$ is differentiable at time $t$, and function $V(\bX(t))$ is differentiable at time $t$. 

\begin{itemize}
\item Case 1: Suppose that $f_1(t) > 0$ and $f_2(t) > 0$. The fluid model equation in \eqref{eq:fluid} can be written in vector form as
\begin{align}\label{eq:fluid2}
\bX(t) = \bX(0) + \blambda t  + \bD^T \bT(t)
\end{align} 
Then,
$$
\dot{f}_1(t) = {\mathbf{\tilde{e}}_1}^T \bM^{-1} \bQ (\blambda + \bD^T \dot{\bT}(t)) = - {\mathbf{\tilde{e}}_1}^T \bD^{-T} \blambda - 1 \leq -\epsilon,
$$
by the utilization condition and the fact that $\dot{T}_1 + \dot{T}_2 = 1$ since LQ scheduling is work-conserving.
Similarly, 
$$
\dot{f}_2(t) = - {\mathbf{\tilde{e}}_2}^T \bD^{-T} \blambda - 1 \leq -\epsilon.
$$
Thus, $\dot{V}(t)\leq -\epsilon$ in this case. 
\item Case 2: Suppose that group 1 is empty, and group 2 is non-empty. That is, $X_1(t) = X_2(t) = 0$ and $\bX(t) \neq \mathbf{0}$. Then, similar to case 1, one can show that $\dot{f}_2(t) \leq -\epsilon$. However, $\dot{f}_1(t)$ is not necessarily strictly negative. Thus, we need to show that at all regular points $t$ in case 2, either $\dot{f}_1(t) < 0$ or $f_2(t) \geq f_1(t)$ (or both). Since $X_1(t) = X_2(t) = 0$, one can write $f_1(t) = c_1 X_3(t) + c_2 X_4(t)$ and $f_2(t) = c_3 X_3(t) + c_4 X_4(t)$ for some positive constants $c_i, ~i=1,2,3,4$. We re-write $f_1$ in the following two expressions.
\begin{align} \label{f1-1}
f_1 = \frac{c_1}{c_3}f_2 + (c_2 - \frac{c_1 c_4}{c_3}) X_4 \\ \label{f1-2}
f_1 = \frac{c_2}{c_4}f_2 + (c_1 - \frac{c_2 c_3}{c_4}) X_3.
\end{align}
Recall that under LQ scheduling, $\dot{T}_3 = 1$ and $\dot{T}_4 = 0$ if $X_3 > X_4$, which implies that $\dot{X}_3 < 0$ and $\dot{X}_4 \geq 0$ in case 2, since $\dot{f}_2 = c_3 \dot{X}_3 + c_4 \dot{X}_4 < 0$. Similarly, by LQ scheduling $X_4 > X_3$ implies that $\dot{T}_4 = 1$, $\dot{T}_3 = 0$, $\dot{X}_4 < 0$, and $\dot{X}_3 \geq 0$ in case 2. Now suppose that $c_1 c_4 < c_2 c_3$. Then, we expand $f_1$ according to \eqref{f1-1}. Then, $X_4 > X_3$ implies that $\dot{f}_1 < -\frac{c_1}{c_3} \epsilon$. We show that if $X_3 \geq X_4$, $f_2 \geq f_1$ in this case as follows. Fixing $X_4 \geq 0$, the function
$$
\frac{f_2}{f_1} = \frac{c_3 X_3 + c_4 X_4}{c_1 X_3 + c_2 X_4} 
$$
is increasing in $X_3$ since $c_1 c_4 < c_2 c_3$ by assumption. Thus, 
$$
\inf_{X_3: ~X_4 \leq X_3} \frac{f_2}{f_1}
$$
is achieved when $X_3 = X_4$. Therefore, it is enough to show that $f_2 \geq f_1$ when $X_3 = X_4$. When $X_1 = X_2 = 0$ and $X_3 = X_4$, we have
\begin{align}\label{x3x4}
\frac{f_2}{f_1} = \frac{\beta[\mu_3^{-1}(q_{33}+ q_{34}) + \mu_4^{-1}(q_{43}+ q_{44})]}{\mu_1^{-1}(q_{13}+ q_{14}) + \mu_2^{-1}(q_{23}+ q_{24})} = \frac{\beta}{a} \geq 1.
\end{align}
In the above derivation, we used \eqref{a} and the fact that $\beta \in [a,b]$.

Similarly, suppose that $c_1 c_4 > c_2 c_3$. Then, we expand $f_1$ according to \eqref{f1-2}. Then, $X_3 > X_4$ implies that $\dot{f}_1 < -\frac{c_2}{c_4} \epsilon$. We show that if $X_4 \geq X_3$, $f_2 \geq f_1$ in this case as follows. Fixing $X_3 \geq 0$, the function 
$$
\frac{f_2}{f_1} = \frac{c_3 X_3 + c_4 X_4}{c_1 X_3 + c_2 X_4} 
$$
is increasing in $X_4$ since $c_1 c_4 > c_2 c_3$ by assumption. Then, 
$$
\inf_{X_4: ~X_3 \leq X_4} \frac{f_2}{f_1}
$$
is achieved when $X_4 = X_3$. Therefore, it is enough to show that $f_2 \geq f_1$ when $X_3 = X_4$, which is already shown in \eqref{x3x4}.

Finally, supposing that $c_1 c_4 = c_2 c_3$, $f_1 / f_2 = c_1 / c_3$. Thus, 
$$\dot{f}_1 = \frac{c_1}{c_3} \dot{f}_2 \leq - \frac{c_1}{c_3} \epsilon.$$

\item Case 3: Suppose that group 2 is empty, and group 1 is non-empty. Then, $\dot{f}_1(t) \leq -\epsilon$. Similar to case 2, we need to show that either $\dot{f}_2(t) < 0$ or $f_1(t) \geq f_2(t)$ (or both). Since in case 3, $X_3 = X_4 = 0$, one can write $f_1 = c'_1 X_1 + c'_2 X_2$ and $f_2 = c'_3 X_1 + c'_4 X_2$ for some positive constants $c'_i, ~i=1,2,3,4$. We have
\begin{align} \label{f2-1}
f_2 = \frac{c'_3}{c'_1}f_1 + (c'_4 - \frac{c'_3 c'_2}{c'_1}) X_2 \\ \label{f2-2}
f_2 = \frac{c'_4}{c'_2}f_1 + (c'_3 - \frac{c'_4 c'_1}{c'_2}) X_1.
\end{align}
Similar to case 2, by LQ scheduling, $X_1 > X_2$ implies that $\dot{X}_1 < 0$ and $\dot{X}_2 \geq 0$. Also, $X_2 > X_1$ implies that $\dot{X}_2 < 0$ and $\dot{X}_1 \geq 0$. Now suppose that $c'_1 c'_4 > c'_2 c'_3$. Then, $X_2 > X_1$ implies that $\dot{f}_2 < -\frac{c'_3}{c'_1} \epsilon$ due to \eqref{f2-1}. We show that if $X_1 \geq X_2$, $f_1 \geq f_2$. Fixing $X_2 \geq 0$, the function
$$
\frac{f_1}{f_2} = \frac{c'_1 X_1 + c'_2 X_2}{c'_3 X_1 + c'_4 X_2} 
$$
is increasing in $X_1$ since $c'_1 c'_4 > c'_2 c'_3$. Thus, 
$$
\inf_{X_1: X_1 \geq X_2} \frac{f_1}{f_2}
$$
is achieved when $X_1 = X_2$. Therefore, it is enough to show that $f_1 \geq f_2$ when $X_1 = X_2$. When $X_3 = X_4 = 0$ and $X_1 = X_2$, we have
\begin{align}\label{x1x2}
\frac{f_1}{f_2} = \frac{\mu_1^{-1}(q_{11}+ q_{12}) + \mu_2^{-1}(q_{21}+ q_{22})}{\beta[\mu_3^{-1}(q_{33}+ q_{34}) + \mu_4^{-1}(q_{41}+ q_{42})]} = \frac{b}{\beta} \geq 1.
\end{align}
In the above derivation, we used \eqref{b} and the fact that $\beta \in [a,b]$.

Similarly, suppose that $c'_1 c'_4 < c'_2 c'_3$. Then, $X_1 > X_2$ implies that $\dot{f}_2 < -\frac{c'_4}{c'_2} \epsilon$ due to \eqref{f2-2}. We show that if $X_2 \geq X_1$, $f_1 \geq f_2$. Fixing $X_1 \geq 0$, the function 
$
\frac{f_1}{f_2}
$
is increasing in $X_2$ since $c'_2 c'_3 > c'_1 c'_4$. Then, 
$$
\inf_{X_2:~X_1 \leq X_2} \frac{f_1}{f_2}
$$
is achieved when $X_2 = X_1$. Therefore, it is enough to show that $f_1 \geq f_2$ when $X_1 = X_2$, which is already shown in \eqref{x1x2}.

Finally, supposing that $c'_1 c'_4 = c'_2 c'_3$, $f_1 / f_2 = c'_1 / c'_3$. Thus, 
$$\dot{f}_2 = \frac{c'_3}{c'_1} \dot{f}_1 \leq - \frac{c'_3}{c'_1} \epsilon.$$
\end{itemize} 
Now taking 
$$\delta = \min(\epsilon,\frac{c_1}{c_3} \epsilon,\frac{c_2}{c_4} \epsilon, \frac{c'_3}{c'_1} \epsilon,\frac{c'_4}{c'_2} \epsilon) > 0,$$
we have $\dot{V}(t) \leq -\delta$. This completes the proof of Theorem \ref{thm:1}.
 
\begin{figure}
\centering
    \includegraphics[width= 0.4 \textwidth]{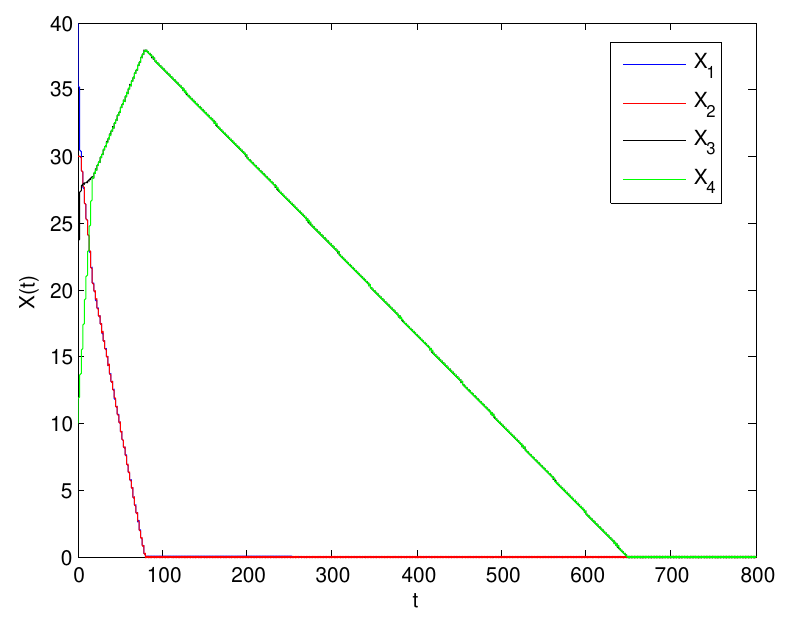}
  \caption{\label{fig:sim1}\textbf{Simulation result.} The figure shows the trajectories of the fluid levels in the queues for Lu-Kumar network under LQ scheduling.} 
\end{figure}

We show the simulation result for the Lu-Kumar network (Figure \ref{fig0}) in the fluid model. In the simulation, service rates are $\bmu = [3,1,1,1]^T$, and arrival rate is $\lambda = 0.4$. Initial fluid levels are $\bX(0) = [40,30,20,10]^T$. Figure \ref{fig:sim1} shows the trajectories of the fluid in different queues versus time. As we can see, 
the fluid level at all queues becomes zero after some finite time.

\section{LDQ Scheduling}\label{sec:LDQ}

\subsection{Policy}

Let $\bar{S}$ be the set of global maxima in the network: 
$$
\bar{S} = \{i: \bar{X}_i = \max_k \bar{X}_k \}.
$$
We define a queue to be ``dominating" if it belongs to the set of global maxima in the network, or if it does not feed another globally maximum-length queue in the network. More formally, for each group $j$, the dominating set of queues $\bar{\calD}_j$ is defined to be 
$$\bar{\calD}_j = \{i \in \calG_j: r_{is}=0 ~ \text{if} ~  \bar{X}_i < \bar{X}_s, \forall s \in \bar{S}\}.$$

\begin{figure}
\centering
    \includegraphics[width= 0.4 \textwidth]{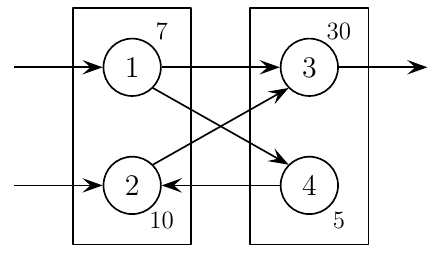}
  \caption{\label{fig5}Illustration of LDQ scheduling.} 
\end{figure}

The scheduling policy is to serve the longest queue in $\bar{\calD}_j$. If $\bar{\calD}_j = \emptyset$, do not serve any queues from group $j$. As an example, see Figure \ref{fig5}, where queue $3$ is the maximum with length $30$. Since jobs leaving both queues $1$ and $2$ can be destined to queue $3$, $\bar{\calD}_1 = \emptyset$ and $\bar{\calD}_2 = \{3,4\}$, so no queues in group $1$ is served, and queue $3$ in group $2$ is served. As we can see, while the policy is not work-conserving, it always serves the queue with maximum-length in the network. The policy is robust to the knowledge of service rates and exact values of routing probabilities. The draw-back in comparison with LQ scheduling is that it needs global knowledge of queue lengths and the topology of the network. 

\subsection{Main result for LDQ scheduling} 
We analyze the fluid model of the system under LDQ scheduling policy. Let $S$ be the set of queues in the network with maximum fluid level: 
$$
S = \{i: X_i = \max_k X_k \}.
$$
Let 
$$
{\calD}_j = \{i \in \calG_j: r_{is}=0 ~ \text{if} ~  X_i < X_s, \forall s \in S\}
$$
be the set of dominating queues in group $j$ in the fluid model. 
Let 
$$
\tilde{S}_j = \{i: i \in {\calD}_j, X_i = \max_{k \in \calD_j} X_k \}
$$ 
be the set of longest dominating queues in group $j$. Define a subset of groups $\calJ \subseteq \{1,2,\ldots,J\}$ as follows:
$$
\calJ = \{j | \exists i \in {\calG}_j \cap S ~ \text{such that} ~ r_{is} = 0, \forall s \in S\}.
$$
In words, $\calJ$ is the set of groups that have a globally maximum-length queue which does not feed another globally maximum-length queue. The scheduling policy implies that  

\begin{align}\label{eq:ldq}
\sum_{i \in \tilde{S}_j} \dot{T}_i(t) &=  1, ~~ \text{if} ~ j \in \calJ, \\ \label{eq:ldq1}
\sum_{i \in \tilde{S}_j} \dot{T}_i(t) & \leq 1,  ~~ \text{if} ~ j \notin \calJ, \\ \label{eq:ldq2}
\sum_{i \in \calG_j \setminus \tilde{S}_j} \dot{T}_i(t) & = 0.
\end{align} 

The fluid model equations for this policy are the one stated in \eqref{eq:fluid} together with \eqref{eq:ldq}--\eqref{eq:ldq2}. 

\begin{thm}\label{thm:2}
A multiclass queueing networks is stable under LDQ scheduling if the utilization condition stated in \eqref{eq:stability} holds and the network is acyclic.
\end{thm}

We prove Theorem \ref{thm:2} using the Lyapunov function $V(\bX)=\max_i \{X_i \}$. Let $|S(t)| = L'(t)$. Let
$S'_j = \calG_j \cap S$ 
be the set of maximum-length (maximum fluid level) queues in group $j$. At time $t$, if $L'=1$, that is the queue with maximum fluid level is unique, $\dot{V}(\bX(t))$ is clearly negative. Here is a proof. Let $S = \{ i \}$. Then, $\dot{T}_i = 1$. Moreover, no jobs will be scheduled destined to the queue with maximum length in the network. Thus, $\dot{V}(\bX(t)) = \lambda_i - \mu_i < 0$. If the maximum is not unique, at a regular point $t$, the drifts of the queues in the set of maxima are equal. (See Lemma 2.8.6. in \cite{dai99}.) That is, $\dot{X}_i(t) = \dot{X}_k(t)$ if $i,k \in S(t)$.

By LDQ scheduling, there is no internal flow to the sub-network of queues in $S$ coming from other queues. That is, if $i \in S$, $\sum_{k \notin S} \mu_k r_{ki} \dot{T}_k = 0$. Thus, we can only consider the sub-network of queues constructed by $S$ to analyze the drift of the Lyapunov function. Let the corresponding drift matrix and exogenous arrival vector to set $S$ be $\bD_{L'}$ and $\blambda_{L'}$. Suppose that $J'$ groups ($J' \leq J$) have queues in $S$, so $S = \cup_{j=1}^{J'} S'_j$ (with some abuse of notation due to relabelling the $J'$ groups by 1 to $J'$). Let $L'_j = |S'_j|$. At a regular point $t$, the drifts of the queues in $S(t)$ are all equal to $\dot{V}(t)$.

\begin{lem}
If $\bX(t) \neq \mathbf{0}$, $\dot{V}(t) < 0$.
\end{lem}

\begin{proof}
Define matrix $\bE_{L' \times J'}$ as the following.
\begin{align}
\bE = [\be_1,\be_2,\dots,\be_{J'}],
\end{align}
where $\be_j$ is a column vector of length $L'$: 
\begin{align}
\be_j = [\mathbf{0}_{\sum_{k=1}^{j-1} L'_k}, \mathbf{1}_{L'_j}, \mathbf{0}_{\sum_{k=j+1}^{J'} L'_k}]^T.
\end{align}

Re-writing (the derivative of) \eqref{eq:fluid} in matrix form together with the fact that at a regular point $t$, if $i \in S$, $\dot{V}(t) = \dot{X}_i(t)$, we have
\begin{equation}\label{eq:matrix2}
\left [ \begin{array}{cc}
- (\bD_{L'})^T & \bE \\
\bE^T & \mathbf{0}_{J' \times J'} 
\end{array}
\right] 
\left [ \begin{array}{c}
\dot{\bT}_{L'} \\
\dot{V} \mathbf{1}_{J'} 
\end{array}
\right] = 
\left [ \begin{array}{c}
\blambda_{L'} \\
\bt_{J'}, 
\end{array}
\right],  
\end{equation}
where vector $\bt_{J'} = [t_j], ~ t_j = \sum_{i \in S'_j} \dot{T}_i \leq 1$. By \eqref{eq:ldq}, at least one of the components of $\bt_{J'}$ is equal to $1$, if the network topology is acyclic. The reason is that the sub-network certainly has an ``output'' queue, which does not feed any other queues in the sub-network. Note that if there is a cycle in the network, such a queue may not exist (See Figure \ref{fig6}). The group which contains the ``output" queue, let's say group $j^*$, has the property that $j^* \in \calJ$; thus, $\sum_{i \in S'_{j^*}} \dot{T}_i = 1$.

\begin{figure}
\centering
    \includegraphics[width= 0.5 \textwidth]{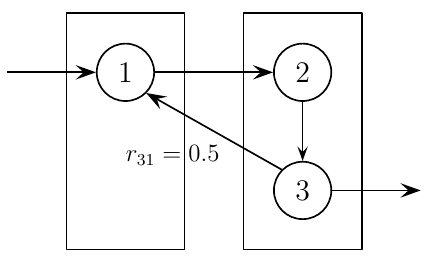}
  \caption{\label{fig6}\textbf{A simple cyclic network.} The figure illustrates that in cyclic networks, there may be no groups in $\calJ$. Note that in this example, if all the queues have equal fluid level, they all belong to the set of global maxima in the network. Now since the network is cyclic, there is no output queue that does not feed any other maximum-length queue.} 
\end{figure}

Solving \eqref{eq:matrix2} using block inverse formula, we have 
\begin{equation}
\dot{V} \mathbf{1}_{J'} = (\bE^T {\bD_{L'}}^{-T} \bE)^{-1} \left(\bE^T   {\bD_{L'}}^{-T} \blambda_{L'}  + \bt_{J'}  \right).
\end{equation}

The vector $ \left(\bE^T {\bD_{L'}}^{-T} \blambda_{L'} + \bt_{J'}  \right)$ is strictly positive in at least one component corresponding to group $j^*$ by the utilization condition of this sub-network (which is a weaker condition than the utilization condition of the whole network).
Since $\bE^T {\bD_{L'}}^{-T} \bE $ is a negative matrix and $\bE^T {\bD_{L'}}^{-T} \bE (\dot{V} \mathbf{1}_{J'})$ is not a negative vector, $\dot{V}$ cannot be non-negative.
\end{proof}

Consequently, $V(\bX(t))$ is a Lyapunov function and the proof of Theorem \ref{thm:2} is complete. 

\begin{figure}[h]
\centering
    \includegraphics[width= 0.4 \textwidth]{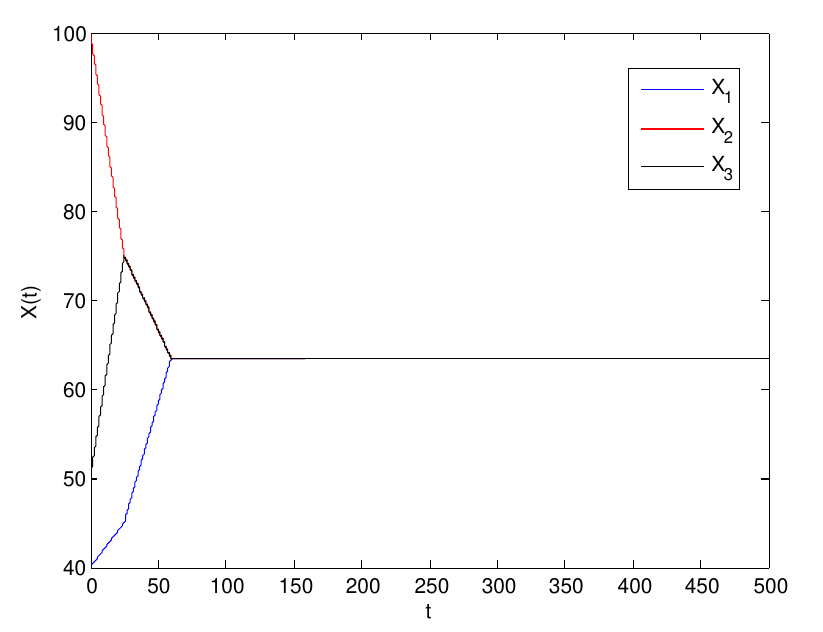}
  \caption{\label{fig7}\textbf{Simulation result.} The figure shows the trajectories of the fluid levels in the queues for the cyclic network shown in Figure \ref{fig6} under LDQ scheduling.} 
\end{figure}

The simulation results show that the network shown in Figure \ref{fig6} is indeed unstable under LDQ scheduling. In this simulation all the service rates are $1$, and the arrival rate is $0.2$. The result of the simulation is shown in Figure \ref{fig7}. As one can see, the fluid levels of the queues do not reach zero. 

\begin{figure}[h]
\centering
    \includegraphics[width= 0.45 \textwidth]{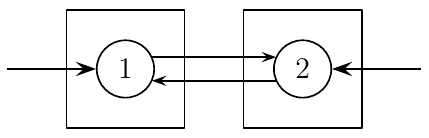}
  \caption{\label{fig8}A simple cycle.} 
\end{figure}

The basic intuition behind this fact is the following. Consider the simple case of Figure \ref{fig8}. We have one queue in each group so no scheduling is needed. However, in the stochastic system, LDQ serves only one of the queues at a time if the queue-lengths are not equal. This illustrates that LDQ may not be throughput-optimal if we have a cycle in the network.

\begin{figure}[h]
\centering
    \includegraphics[width= 0.4 \textwidth]{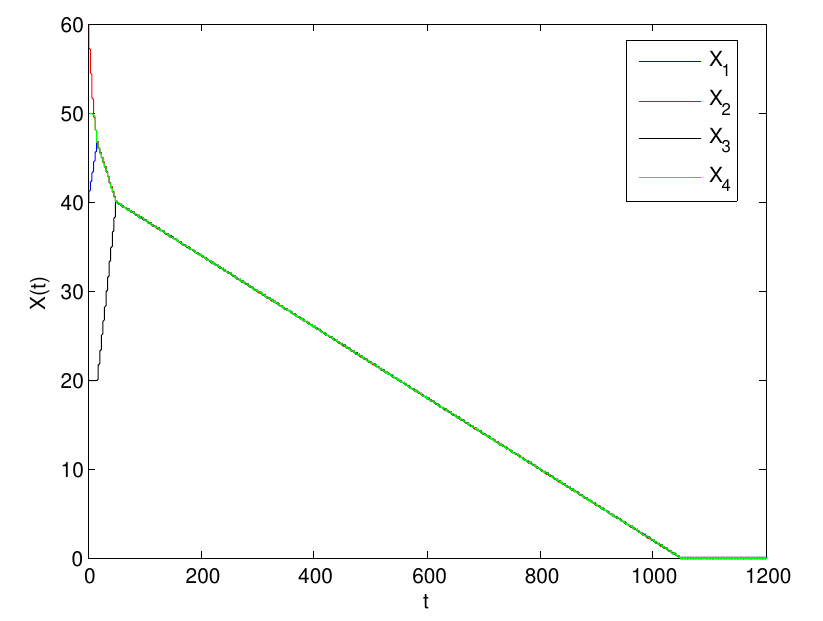}
  \caption{\label{fig9}\textbf{Simulation result.} The figure shows the trajectories of the fluid levels in the queues for Lu-Kumar network under LDQ scheduling.} 
\end{figure}

In the end, we see the simulation result for the Lu-Kumar network under LDQ scheduling. The network topology is acyclic, so we expect LDQ to be stable. The simulation result shown in Figure \ref{fig9} verifies this. In this simulation, all the service rates are $1$, and the arrival rate is $0.4$. 
As one can see, the fluid level at all queues becomes zero after some finite time.

\acks
This work is supported by MURI grant BAA 07-036.18. We thank the anonymous referees for many helpful remarks.

\end{document}